\documentclass[11pt]{{amsart}}

\usepackage[all]{xy}
\usepackage[ansinew]{inputenc}
\usepackage[psamsfonts]{amssymb}
\theoremstyle{plain}
\newtheorem{theorem}{Theorem}
\newtheorem{proposition}{Proposition}
\newtheorem{lemma}[proposition]{Lemma}

\newtheorem{problem}{Open Problem}

\theoremstyle{definition}
\newtheorem{definition}[proposition]{Definition}

\theoremstyle{remark}
\newtheorem{remark}[proposition]{Remark}

\newcommand{\secref}[1]{Section~\ref{#1}}

\newcommand{\thmref}[1]{Theorem~\ref{#1}}
\newcommand{\propref}[1]{Proposition~\ref{#1}}

\newcommand{\defref}[1]{Definition~\ref{#1}}

 \def\bk{{l\!k}} 
\def\Q{\mathbb Q}
\def\F{\mathbb F}
\def\R{\mathbb R}
\def\C{\mathbb C}
\def\Z{\mathbb Z}
\def\cE{\mathcal E}
\def\cF{\mathcal F}
\def\cG{\mathcal G}
\def\cP{\mathcal P}
\title{The cohomology algebra of unordered configuration spaces}
\author{Yves F\'elix}
\address{D\'epartement de Mathematiques\\
        Universit\'e Catholique de Louvain\\
         2, Chemin du Cyclotron\\
          1348 Louvain-La-Neuve\\
         Belgium}
\email{felix@uranus.sc.ucl.ac.be}
\author{Daniel Tanr\'e}
\address{D\'epartement de Mathematiques\\
         UMR 8524\\
         Universit\'e de Lille~1\\
         59655 Villeneuve d'Ascq Cedex\\
         France}
\email{Daniel.Tanre@agat.univ-lille1.fr}

\begin{document}

\begin{abstract}
Given an $N$-dimensional compact manifold $M$ and a field $\bk$, F.~Cohen
and L.~Taylor have constructed a spectral sequence, $\cE(M,n,\bk)$,
converging to the cohomology of the space of ordered configurations of
$n$ points in $M$. The symmetric group $\Sigma_n$ acts on this spectral
sequence giving a spectral sequence of $\Sigma_n$ differential graded
commutative algebras. Here, we provide an explicit description of the
invariants algebra
$(E_1,d_1)^{\Sigma_n}$ of the first term of $\cE(M,n,\Q)$. We apply this
determination in two directions:

-- in the case of a complex projective manifold or of an odd dimensional
manifold $M$, we obtain the cohomology algebra $H^*(C_n(M);\Q)$ of the
space of unordered configurations of $n$ points in $M$ (the concrete
example of $P^2(\C)$ is detailed),

-- we prove the degeneration of the spectral sequence formed of the
$\Sigma_n$-invariants $\cE(M,n,\Q)^{\Sigma_n}$ at level~2, for any
manifold $M$.

These results use a transfer map and are also true with coefficients in a
finite field $\F_p$ with $p>n$.
\end{abstract}

\maketitle
Given an $N$-dimensional compact manifold, $M$, define
the space $F(M,n)$ of  \emph{ordered configurations} of $n$ points in
$M$, as
$$F(M,n) = \,\{\, (x_1, x_2, \ldots , x_n)\, \vert \, x_i \neq
x_j\,, \mbox{for } i\neq j\,\}\,.$$ The symmetric group $\Sigma_n$
acts freely on $F(M,n)$ by permutation of coordinates. Here we are
interested in the  orbit space
$C_n(M) = F(M,n)/
\Sigma_n$ called the space of \emph{unordered configurations} of $n$
points in $M$. The determination of the homotopy type, even of the
rational homotopy type, of the spaces $F(M,n)$ and $C_n(M)$ is a
hard problem.

Let $\bk$ be a field. The principal tool for the determination of the
cohomology algebra of
$F(M,n)$ is a spectral sequence $\cE(M,n,\bk)$ of
graded commutative algebras constructed by  F.~Cohen and
L.~Taylor in
\cite{CT}.  This spectral sequence converges to
$H^*(F(M,n);\bk)$ and has a $E_1$-term that we briefly recall now:
$$ (E_1,d_1) = ((H^*(M;\bk)^{\otimes n} \otimes \land_{1\leq
i,j\leq n} e_{ij})/I, d_1)\,.$$
The elements of $(\otimes^n H(M;\bk))^t$ have bidegree $(t,0)$, the
$e_{ij}$ bidegree $(0,N-1)$ and the differential denoted by $d_1$ is the
first non-zero differential. The symbol
$\land V$ denotes the free commutative graded algebra on the graded vector
space
$V$. The ideal
$I$ is generated by the elements
$e_{ij}e_{jk} + e_{jk}e_{ki} + e_{ki}e_{ij}$, \hskip .2cm
$e_{ij}^2$, \hskip .2cm
$e_{ij} - (-1)^Ne_{ji}$  and
$(p_i^*(a)-p^*_j(a))\otimes e_{ij}$ where $p_i : M^n \to M$ is the
projection on the $i^{th}$ factor. Finally,
$$d_1(e_{ij}) = p_{ij}^*(\Delta),$$
where $p_{ij} : M^n \to M^2$ is the projection on the $i^{th}$ and
the $j^{th}$ factors and $\Delta \in H^N(M^2;\bk)$ is the diagonal
class.

The symmetric group $\Sigma_n$ acts on $(E_1,d_1)$ by permuting the
coordinates in $H^*(M;\bk)^{\otimes n}$ and on
$\land_{1\leq
i,j\leq n} e_{ij}$ by $\sigma(e_{ij})=e_{\sigma(i)\sigma(j)}$. With
this action, the spectral sequence $\cE(M,n,\bk)$ becomes  a spectral
sequence of
$\Sigma_n$-algebras and creates a spectral sequence
$\cE(M,n,\bk)^{\Sigma_n}$ consisting of the
$\Sigma_n$-invariants of each stage.
With the basic results recalled in \secref{sec:equivarianthomology},
observe that, for
$\bk=\Q$, or for $\bk=\F_p$ and $p>n$, the spectral sequence
$\cE(M,n,\bk)^{\Sigma_n}$ converges to
$H^*(C_n(M);\bk)$.

 This work is concerned with the
study of the degeneration of the spectral sequence
$\cE(M,n,\bk)^{\Sigma_n}$
and its
implication in the determination of the cohomology algebra of
$H^*(C_n(M);\bk)$.

First recall  the particular situation of a smooth projective complex
manifold $M$ and $\bk =\Q$.   I. Kriz
(\cite{K}) and B. Totaro (\cite{T})   prove  that the differential
graded algebra
$(E_1,d_1)$ is quasi-isomorphic to the Sullivan minimal model of
the space $F(M,n)$ which  means that $(E_1,d_1)$ contains all the
rational homotopy type of $F(M,n)$.  In
particular, in this case, the spectral sequences $\cE(M,n,\Q)$
and $\cE(M,n,\Q)^{\Sigma_n}$ collapse at level~2 and
$H^*(C_n(M);\Q)\cong H(E_1,d_1)^{\Sigma_n}$ as algebra
(see \cite[Theorem 8 and Corollary 8c]{FM}, \cite[Remark 1.3]{K} and
\cite[Corollary~1 of Theorem~5]{T}).
For more general manifold we
do not have a model for the rational homotopy type except for $n =
2$ by a work of P. Lambrechts and D. Stanley (\cite{LS}).

As we see in the particular case of a projective complex manifold and
$\bk=\Q$, we have first to understand the invariants algebra
$(E_1,d_1)^{\Sigma_n}$. We do that for any field $\bk$ and prove in
\thmref{thm:invariantsk} of \secref{sec:invariantCohenTaylor} that, for
$N$ even, we have an isomorphism
$$E_1^{\Sigma_n} \cong \Gamma_n^H\,,\,{\text with}\;
\Gamma_n^H =
\oplus_{r=0}^{[n/2]}
\Gamma^{n-2r}(H) \otimes \Gamma^{r}(s^{N-1}H),$$
where $\Gamma^j(V)\subset T^j(V)$ denotes the invariants subspace
for the action of the symmetric group $\Sigma_j$ by permutation of the
factors. The law structure induced on $\Gamma_n^H$ by this isomorphism
is described in
\thmref{thm:invariantskalg}. In the case $\bk=\Q$ or
$\bk=\F_p$ with $p>n$, this result becomes

\begin{theorem}\label{thm:invariantQ}
Let $M$ be a
$N$-dimensional compact manifold with $N$ even and let $\bk=\Q$ or
$\bk=\F_p$ with
$p>n$. We denote by
$\bullet$ the  multiplication law of the cohomology algebra
$H = H^*(M;\Q)$. Then there is an isomorphism
of graded algebras
$$(E_1,d_1)^{\Sigma_n} \cong (C_n^H,d)\,,\,{\text with}\; C_n^H =
\left(
\oplus_{r=0}^{[n/2]}
\land^{n-2r}(H) \otimes \land^{r}(s^{N-1}H),d\right).$$
where the multiplication law on $C_n^H$ depends only on the
multiplication law $\bullet$ on $H$ and is described in
\thmref{thm:multiplicationofcnH}.
The differential
$$d\colon \land^{n-2r}(H) \otimes \land^{r}(s^{N-1}H)\rightarrow
\land^{n-2r+2}(H) \otimes \land^{r-1}(s^{N-1}H)$$
 is determined by
  $$d ( e_0 \land \ldots\land e_0 \otimes  s^{N-1} a) =
  \frac{1}{2}\sum_k (-1)^{\vert b_k'\vert }
  e_0 \land\ldots\land  e_0 \land(a\bullet b_k) \land b_k'\,,$$
where $e_0$ is the unit of $H$, the    $\{ b_i\}$ form an homogeneous
basis for $H$ and the
$\{ b_i'\}
 $ are the dual basis for Poincar\'e duality.
\end{theorem}

The differential $d$ of a general element is expressed in
\defref{def:generald}. If we consider
$C^H=\oplus_n C_n^H= \land H\otimes \land (s^{N-1}H)$, it appears from
\defref{def:generald} that the differential $d$ is also a differential of
algebra for the \emph{usual product} on $\land H\otimes \land
(s^{N-1}H)$. Denote by $\Gamma_{sec}(M)$ the space of sections of the
sphere bundle associated to the tangent bundle of $M$. The commutative
differential graded algebra $\land H\otimes \land
(s^{N-1}H)$ with the usual product is the model of $\Gamma_{sec}(M)$
already determined in \cite{FT}. The fact that our $C_n^H$ fit together
and give this model corresponds to the link between
$\Gamma_{sec}(M)$
and $\vee_{n\geq 0} C_n(M)$ revealed in \cite{McD}, see also \cite{B}.

\emph{In the case of a complex projective manifold $M$,}
the Kriz-Totaro result recalled above gives an isomorphism of algebras
between
$H^*(C_n(M);\Q)$ and $H(C_n^H,d)$. This description is
effective and in
\secref{sec:configCP2}, we can determine explicitely the rational
cohomology algebra of the unordered configuration spaces of $n$ points
in $P^2(\C)$:

\begin{theorem}\label{thm:configCP2}
The rational cohomology of the
 unordered configuration spaces of $n$ points into $P^2(\mathbb
 C)$ is given by
 $$
 \renewcommand{\arraystretch}{1.6}
 \begin{array}{lll}
 n=1,2 \hspace{5mm}\mbox{}& \land x/x^3 & \vert x\vert = 2\,;\\
 n=3 & \land (x,y)/ (x^3, yx^2) & \vert x\vert = 2\,, \vert y\vert
 = 3\,;\\
 n\geq 4 & \land (x,y)/x^3 & \vert x\vert = 2\,, \vert y\vert
 = 3\,.
 \end{array}
\renewcommand{\arraystretch}{1}
 $$
\end{theorem}

In the more general case of a real manifold, recall that the
homology vector space $H_*(C_n(M);\bk)$ has been described, when $N$ is
odd,  by C.-F.~B\"odigheimer, R.~Cohen and L.~Taylor in \cite{BCT},
and, when $N$ is even and $\bk=\Q$, by Y.~F\'elix and J.-C.~Thomas in
\cite{FT}. From these results, we deduce:

\begin{theorem}\label{thm:collapseQ}
For any manifold $M$, the spectral sequence formed of the
$\Sigma_n$-invariants,
$\cE(M,n,\Q)^{\Sigma_n}$, collapses at level~2.
\end{theorem}

Observe that this theorem contrasts with a recent result of the first
author and J.-C.~Thomas. Answering a question of M.~Bendersky and
S.~Gitler (\cite{BG}), they prove that there exist manifolds $M$ for which
the spectral sequence
$\cE(M,n,\Q)$ does not collapse at level~2, see \cite{FT2}.

\smallskip
The result of \cite{BCT} being true over the field $\F_p$, we have also

\begin{theorem}\label{thm:casimpair}
Let $M$ be an odd dimensional manifold
and let $\bk=\Q$ or
$\bk=\F_p$ with
$p>n$. The spectral
sequence $\cE(M,n,\bk)^{\Sigma_n}$ collapses at level~2 and the
cohomology algebra of the unordered configuration space is isomorphic to
$\land^n(H)$.
\end{theorem}

As a conclusion, observe that two cases are not covered by our study
and they lead to

\begin{problem} Does the spectral sequence
$\cE(M,n,\F_p)^{\Sigma_n}$ collapse for an even dimensional manifold?
\end{problem}

What is missing for answering is the knowledge of the homology
$H_*(C_n(M);\F_p)$ as we have it over $\Q$ in the even dimensional case.
Notice now that
\thmref{thm:collapseQ} gives an information on the algebra structure of
$H^*(C_n(M);\Q)$ only up to a filtration.

\begin{problem}
Can we recover the
algebra structure of
$H^*(C_n(M);\Q)$ in all cases?
\end{problem}

We know the algebra structure when $M$ is odd dimensional,
when $M$ is projective complex or when $n=2$ (see \cite{LS}) but not in
the general case.

\smallskip
In \secref{sec:equivarianthomology}, we recall some basic facts
concerning equivariant homology and cohomology.
\secref{sec:invariantofconfigurationinRN} and
\secref{sec:invariantCohenTaylor} are concerned with the space of
$\Sigma_n$-invariants of $H^{(n-1)(N-1)}(F({\mathbb
R}^N,n);\bk)$  and $E_1$ respectively.
\secref{sec:CnH} is devoted to the differential graded
algebra $(C_n^H,d)$ and the proof of \thmref{thm:invariantQ}. It contains
also the proof of
\thmref{thm:collapseQ}
in the even dimensional case.
As an example we give in
\secref{sec:configCP2} a complete description of the cohomology algebra of
the configuration spaces of $n$ points in the complex projective space
$P^2(\mathbb C)$. In \secref{sec:oddcase}, we study the case of odd
dimensional manifolds and prove \thmref{thm:casimpair}.

\section{Equivariant homology}\label{sec:equivarianthomology}

Let $G$ be a finite group and let $\bk$ be a field. A
\emph{$G$-cdga} is a commutative differential graded algebra
$(A,d_A)$ with $H^0(A,d_A) = \bk$, on which $G$ acts by cdga
maps.
The invariant subspace $(A,d_A)^G$ defines a sub cdga of $(A,d_A)$.
Recall some basic facts on these objects:

\begin{proposition}\label{prop:gquasiiso}
Let $\bk$ be a field of characteristic that does not divide the order of
$G$ and let $f : (A,d_A) \to
(B,d_B)$ be a $G$-equivariant quasi-isomorphism. Then, we have

1) $f^G :
(A,d_A)^G \to (B,d_B)^G$ is also a quasi-isomorphism;

2) $H\left((A,d_A)^G\right)= \left(H(A,d_A)\right)^G$.
\end{proposition}

\begin{proof} 1) Let $a\in (A,d_A)^G$ be a
cocycle and suppose that $f^G(a) = d(b)$, then $f(a) = db$ and
$a=d(c)$. Since $a$ is invariant, $a= d\left(\frac{1}{\vert
G\vert} \sum_{g\in G} g\cdot c\right)$. This shows that $H(f^G)$ is
injective.

Let now $b\in  B^G$ be a cocycle. Since $f$ is a
quasi-isomorphism, there is a cocycle $a \in A$ and an element
$c\in B$ such that $f(a) = b+d(c)$. Therefore $f\left(\frac{1}{\vert
G\vert } \sum_{g\in G}g\cdot a\right) = b + d\left(\frac{1}{\vert G\vert}
\sum_{g\in G} g\cdot c\right)$. Thus $H(f^G)$ is also surjective.

2) The second statement  is a consequence of the existence of
the canonical inclusion $A^G\rightarrow A$ and of $\sigma\colon
A\rightarrow A^G$,
$\sigma(a)= \sum_{g\in G} g \cdot a$.
\end{proof}

Let $X$ be a simplicial complex with a (simplicial) action of $G$. Recall
from G.~Bredon (\cite[Page~115]{Br}) that the action is \emph{regular} if,
for any
$g_0,\ldots,g_n$ in $G$ and simplices
$(v_0,\ldots,v_n)$, $(g_0v_0,\ldots,g_nv_n)$
of $X$, there exists an element $g\in G$ such that $gv_i=g_iv_i$ for all
$i$. By \cite[Proposition~1.1, Page 116]{Br}, the induced action on the
second barycentric subdivision is always regular.
Here, \emph{by a $G$-space, we mean a connected simplicial complex on
which
$G$ acts regularly.}

\smallskip
Denote by $C(X)$ the oriented chain complex of $X$ and observe that
$C(X)$ is a module over the group ring $\Z[G]$ of $G$. The canonical
simplicial map
$\rho\colon X\to X/G$ induces $\rho_{\ast}\colon C(X)\to C(X/G)$. Define
now
$\sigma\colon C(X)\to C(X)$, $c\mapsto \sum_{g\in G} g c$. One has
${\text Ker} \,\sigma = {\text Ker} \,\rho_*$. Therefore $\sigma$ induces
$\overline{\sigma}$
$$\xymatrix{
C(X)\ar[r]^-{\sigma}\ar[d]_-{\rho_*}&C(X)\\
C(X/G)\ar[ru]_-{\overline{\sigma}}&\\
}$$
such that $\overline{\sigma}\circ \rho_*=\sigma$. One can prove:

\begin{proposition}[{\cite[Page 120]{Br}}]
Let $\bk$ be a field of characteristic that does not divide the order of
$G$ and let $X$ be a $G$-space. Then there are isomorphisms
$$C_*(X/G;\bk)\cong C_*(X;\bk)^G
\;{\text and}\;
C^*(X/G;\bk)\cong C^*(X;\bk)^G\,.$$
\end{proposition}

In the case $\bk=\Q$, when $X$ is a $G$-space the group
$G$ acts on the Sullivan algebra of PL-forms on $X$, $A_{PL}(X)$,
(\cite{S}). A similar argument identifies
$A_{PL}(X/G)$ with
$A_{PL}(X)^G$.
Moreover, a $G$-cdga $(A,d_A)$ admits a minimal model
$$\varphi : (\land V,d) \stackrel{\simeq}{\longrightarrow}
(A,d_A)$$ with an action of  $G$ on $(\land V,d)$ making  $\varphi$
$G$-equivariant (\cite{GHV}). This model  is unique
up to
$G$-isomorphisms. We call it \emph{a $G$-minimal model} and
\propref{prop:gquasiiso} implies

\begin{proposition}Let $(\land V,d)
\stackrel{\cong}{\to} A_{PL}(X)$ be a $G$-minimal model of the
$G$-space $X$,  then the fixed point set $(\land V ,d)^G$
is a model for $A_{PL}(X/G)$.
\end{proposition}

This means that the cdga's $(\land V ,d)^G$ and
$A_{PL}(X/G)$ have the same minimal model.

\section{$\Sigma_n$-invariants of $H^{(n-1)(N-1)}(F({\mathbb
R}^N,n);\bk)$, $N$ even}\label{sec:invariantofconfigurationinRN}

The computation of $H^*(F({\mathbb R}^N,n);\bk)$ has been realized
by F.~Cohen in \cite{C},
$$H^*(F({\mathbb R}^N,n);\bk) = \land_{1\leq i,j\leq n}
e_{i,j} /I$$ where $I$ is the ideal generated by the elements
$e_{i,j}-(-1)^N e_{j,i}$ and the elements $e_{i,j}e_{j,k} +
e_{j,k}e_{k,i} + e_{k,i}e_{i,j}$. Here the elements $e_{i,j}$ have
all degree $N-1$.  A basis of the cohomology is given by the
products $e_{i_1,j_1}e_{i_2,j_2}\ldots e_{i_r,j_r}$ where for each
$s$, $i_s<j_s$ and $j_1<j_2< \ldots <j_r$.
The permutation group $\Sigma_n$ acts naturally on $\land
(e_{i,j})$ by $\sigma (e_{i,j})  = e_{\sigma (i), \sigma (j)}$.

\smallskip
In
this section, we prove

\begin{theorem}\label{thm:invariantsRN}
For $n\geq 3$, one has $H^{(n-1)(N-1)}(F({\R}^N,n);\bk)^{\Sigma_n}=0$.
\end{theorem}

In \cite{O}, E. Ossa makes a decomposition of the homology of the
configuration space $F(R^N,n)$ as a module on the ring group
$\Z[\Sigma_n]$, see also \cite[Lemma~5.2]{AJ}. As in \cite{O}, our
starting point is similar: we replace the study of configuration spaces in
the setting of trees.

 We denote by ${\mathcal G}_n$ the set of \emph{connected} trees with
$n$ vertices, denoted $1, 2, \ldots , n$, and with edges $(i,j)$
naturally oriented from $i$ to $j$ if $i<j$. For sake of
simplicity we will always denote the edges    $(i,j)$ with $i<j$.

Define $G_n$ as the the quotient of the
vector space over $\bk$, constructed on the elements of ${\mathcal
G}_n$, by the following relation: if a graph $G$ contains the
edges  $(i,k)$ and $(j,k)$ with $i<j$, then $G = G' -G''$, where
$G'$ is obtained by replacing the above edges by the edges $(i,j)$
and $(j,k)$ and $G''$ is obtained by replacing the above edges by
the edges $(i,k)$ and $(i,j)$, i.e.:
$$\vcenter{\xymatrix{
&k&&&k&&&k\\
j\ar@{-}|@{>}[ur]&&=&j\ar@{-}|@{>}[ur]&&-&j&\\
&i\ar@{-}|@{>}[uu]&&&i\ar@{-}|@{>}[ul]&&&
i\ar@{-}|@{>}[ul]\ar@{-}|@{>}[uu]\\
}}\;.$$
Observe that an iteration of this relation allows a choice of
representing elements such that each vertex is the end of at most one
edge.

The group
$\Sigma_n$ acts on
${\mathcal G}_n$ par permutation of vertices and this action induces an
action on the vector space $G_n$. Since a tree is completely
determined by the sequence of its
 edges, we associate to a tree with  edges
 $(i_s,j_s)$, $s=1,\ldots,n-1$, the element $e_{i_1,j_1}e_{i_2,j_2}\ldots
e_{i_{n-1},j_{n-1}}$ of $\land
 (e_{i,j})$. This defines an isomorphism of $\Sigma_n$-vector spaces
 $$G_n \stackrel{\cong}{\longrightarrow} H^{(n-1)(N-1)} (F({\mathbb
R}^N,n);\bk)\,.$$
The proof of \thmref{thm:invariantsRN} is reduced to the study of
$G_n^{\Sigma_n}$.

\medskip
 We now consider
 the stabilizator $\Sigma_{n-1}$ of $1$ in $\Sigma_n$ and we denote by $C_n$ the graph in
 ${\mathcal G}_n$ whose edges are $(1,2), (1,3), \ldots , (1,n-1)$.

\begin{lemma}  For $n\geq 3$, one has $G_n^{\Sigma_{n-1}} = \bk \; C_n$.
\end{lemma}

\begin{proof} We suppose by induction on $n$
 that this is true for $k\leq n-1$. We then decompose $G_n$ as a direct
 sum
 $$G_n =  \oplus_{p\leq n-1} G_{n,p}\,,$$
 where $G_{n,p}$ is the sub-vector space generated by   the
 trees $T$ whose  components $T_1$ and $T_2$ of $1$ and $2$
 in $T\backslash (1,2)$ contain  respectively $n-p$ and $p$
vertices. The wedge $T_1 \vee T_2$ injects into $T$ and we
identify each $T_i$ with its image.

Let $a\in G_n^{\Sigma_{n-1}}$, then $a=\sum_p a_p$, $a_p\in
G_{n,p}$ and each $a_p$  is invariant by the subgroup
$\Sigma_{n-2}$ that fixes the vertices $1$ and $2$. In particular,
via the injection $T_1 \hookrightarrow T$, the group
$\Sigma_{n-p-1}$ acts on $G_{n,p}$, and   by restriction we have a
surjective map
$$G_{n,p}^{\Sigma_{n-p-1}} \to G_{n-p}^{\Sigma_{n-p-1}}\,.$$
By induction this implies that $a_p$ is a linear combination of
trees in which $T_1 = C_{n-p}$. That means that the trees in this
decomposition have the same form but not the same indexing on vertices.
In the same way, we can suppose that
$T_2 = C_p$. We consider now the index~3 and we can write
$$a = \sum_i \alpha_i a_i + \sum_j \beta_j a_j'$$
where   the $a_i$  are trees $T$ with $T_1 =
C_{n-p}$, $T_2 = C_p$ and $3 \in T_2$
$$
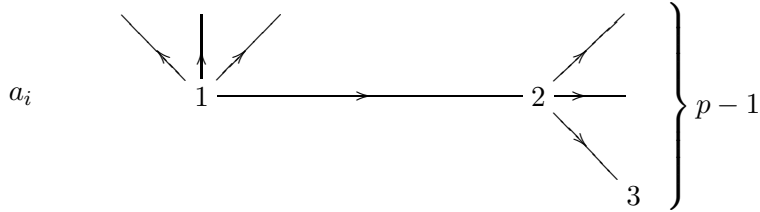
\left.
\vcenter{\xymatrix{
&&& &&& &\\
a_i&&1\ar@{-}|@{>}[ur]\ar@{-}|@{>}[u]\ar@{-}|@{>}[ul]\ar@{-}|@{>}[rrrr]&&&&
2
\ar@{-}|@{>}[ru]\ar@{-}|@{>}[r]\ar@{-}|@{>}[rd]&\\
&&& &&&
&3}}
\;\;\right\} p-1
$$
and the $a_j'$ are trees $T$ with $T_1 =
C_{n-p}$, $T_2 = C_p$ and $3\in T_1$.
$$
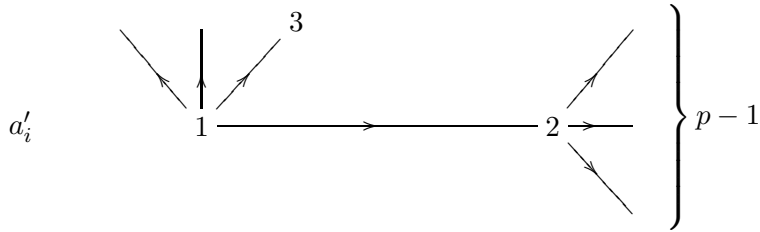
\left.
\vcenter{\xymatrix{
&&&3 &&& &\\
a'_i&&1\ar@{-}|@{>}[ur]\ar@{-}|@{>}[u]\ar@{-}|@{>}[ul]\ar@{-}|@{>}[rrrr]
&&&&
2
\ar@{-}|@{>}[ru]\ar@{-}|@{>}[r]\ar@{-}|@{>}[rd]&\\
&&& &&&
&}}
\;\;\right\} p-1
$$
Denote by $\sigma$ the permutation
$(2,3)$ of the vertices $2$ and $3$. Then the image by $\sigma$ of
the graph $a_i$ is the linear combination $a_{i,1}-a_{i,2}$ with
$$
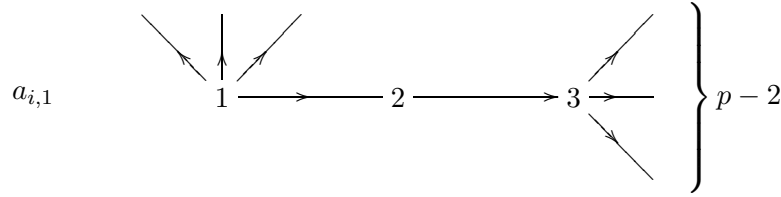
\left.
\vcenter{\xymatrix{
&&& &&& &\\
a_{i,1}&&
1\ar@{-}|@{>}[ur]\ar@{-}|@{>}[u]\ar@{-}|@{>}[ul]\ar@{-}|@{>}[rr]
&&2\ar[rr]&&
3
\ar@{-}|@{>}[ru]\ar@{-}|@{>}[r]\ar@{-}|@{>}[rd]&\\
&&& &&&
&}}
\;\;\right\} p-2
$$
and
$$
\left.
\vcenter
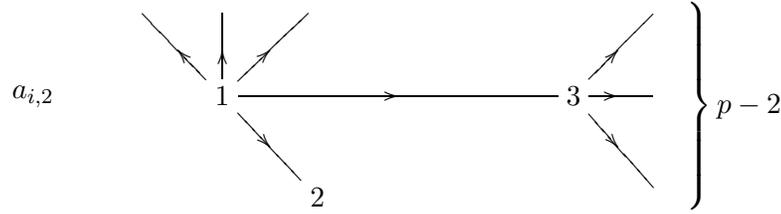
{\xymatrix{
&&& &&& &\\
a_{i,2}&&
1\ar@{-}|@{>}[ur]\ar@{-}|@{>}[u]\ar@{-}|@{>}[ul]\ar@{-}|@{>}[rrrr]
\ar@{-}|@{>}[dr]&&&&
3
\ar@{-}|@{>}[ru]\ar@{-}|@{>}[r]\ar@{-}|@{>}[rd]&\\
&&&2 &&&
&}}
\;\;\right\} p-2
$$
\vspace{1cm} \vspace{2mm} The image of $a_j'$ by $\sigma$ is a graph of
the form
$$
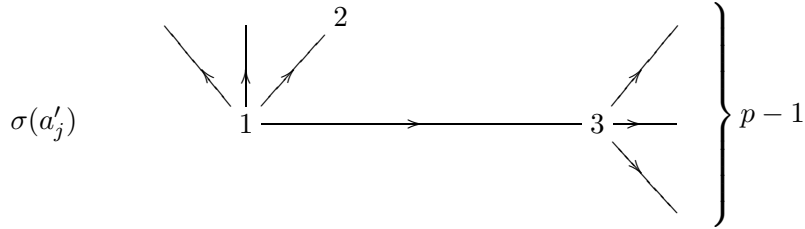
\left.
\vcenter{\xymatrix{
&&&2 &&& &\\
\sigma(a'_j)&&1\ar@{-}|@{>}[ur]\ar@{-}|@{>}[u]\ar@{-}|@{>}[ul]\ar@{-}|@{>}[rrrr]
&&&&
3
\ar@{-}|@{>}[ru]\ar@{-}|@{>}[r]\ar@{-}|@{>}[rd]&\\
&&& &&&
&}}
\;\;\right\} p-1
$$
If we decompose each element $a_{i,1}$, $a_{i,2}$, $\sigma(a'_j)$ in
the direct sum $\oplus \,G_{n,p}$ one has
$a_{i,2}\in G_{n,1}$,
$\sigma (a'_j)\in G_{n,1}$ and
$a_{i,1}\in G_{n,p}$. Therefore, the equality
$a=\sigma(a)$
implies $\alpha_i=0$ for $a_i\in G_{n,p}$, $p\geq 3$. But now, by a
decreasing induction on $p$, one notices $\beta_j=0$ for $a'_j\in
G_{n,p}$,
$p\geq 2$. Thus,
$$\xymatrix@=16pt{
&&&4&5&n&& &&2&3&&\\
a&=&\alpha&1\ar@{-}|@{>}[u]\ar@{-}|@{>}[ur]\ar@{-}|@{>}[urr]
\ar@{-}|@{>}[rr]&&2\ar@{-}|@{>}[r]&3
&+&\beta&1\ar@{-}|@{>}[rr]\ar@{-}|@{>}[u]\ar@{-}|@{>}[ur]&&n
}$$
and
$$\xymatrix@=16pt{
&&&4&5&n&& &&2&3&&\\
\sigma(a)&=&\alpha&1\ar@{-}|@{>}[u]\ar@{-}|@{>}[ur]\ar@{-}|@{>}[urr]
\ar@{-}|@{>}[rr]&&2\ar@{-}|@{>}[r]&3
&+&(\beta-\alpha)&1\ar@{-}|@{>}[rr]\ar@{-}|@{>}[u]\ar@{-}|@{>}[ur]&&n.
}$$
Therefore $\alpha = 0$ and $G_n^{\Sigma_{n-1}} =
\bk \;C_n$.
\end{proof}

\begin{proof}[Proof of \thmref{thm:invariantsRN}] Let $a \in
G_n^{\Sigma_n}$, then $a=\alpha C_n$. Denote $\tau = (1,2)$. Since
$\tau (a) = a$,
$\alpha$ must be equal to $0$ .
\end{proof}

\section{$\Sigma_n$-invariants of the Cohen-Taylor algebra
$E_1$, $N$ even.}\label{sec:invariantCohenTaylor}

Let $M$ be a compact even dimensional manifold of cohomology algebra
$H=H^*(M;\bk)$.
The symmetric group $\Sigma_j$ acts on $H^{\otimes j}$ by permutation of
factors and we denote by
 $\Gamma^j H \subset H^{\otimes j}$ the subspace of
 $\Sigma_j$-invariant elements.

Recall from \secref{sec:invariantofconfigurationinRN},
the existence of a correspondence between connected trees with $n$
vertices and
$H^{(n-1)(N-1)}(F(M,n);\bk)$. This correspondence can be extended to the
$E_1$ term of the Cohen-Taylor spectral sequence. Grants to the relation
$(p_i^*(a)-p^*_j(a))\otimes e_{ij}$, one has
$$E_1= \oplus_{g\in\cF_n} H^{\otimes l(G)} \cdot G,$$
where $\cF_n$ denotes the set of \emph{forests} (i.e. disjoint union of
trees) with $n$ vertices which are the end of at most one edge
and where
$l(G)$ denotes the number of components of $G$. We will make more precise
this isomorphism in a particular case of forests we are interested in.

\smallskip
Let $0\leq r\leq n/2$ be an integer and denote by ${\mathcal P}_r$ the
set of $r$-uples of disjoint pairs $I=\{ (i_1,j_1),
(i_2,j_2),\ldots , (i_r,j_r)\}$ with $i_k<j_k$ for $k\leq r$ and
$1\leq i_1<i_2 < \ldots <i_r<j_r\leq n$. (Disjoint pairs means $i_u\neq
j_v$, $j_u\neq j_v$ for all $u$ and $v$.)

To any $I\in {\mathcal P}_r$, we associate the element
$e_I = e_{i_1,j_1}\cdots e_{i_r,j_r}$ and we note $I_1 = \{
i_1,
\ldots , i_r\}$, $I' = \{ i_1, \ldots , i_r, j_1, \ldots , j_r\}$.
For any fixed $I\in {\mathcal P}_r$, we define a  linear map
$$ \gamma_I : T^r(s^{N-1} H) \to T^r(H) \cdot e_I$$
by
$\gamma_I (s^{N-1}x_1 \otimes \cdots \otimes s^{N-1}x_r) =
(-1)^{(N-1)(\vert x_1\vert + 2\vert x_2\vert + \cdots + r\vert
x_r\vert)} x_1\otimes \cdots \otimes x_r\cdot e_I$. We define a
$\Sigma_r$-action on the domain and the target of $\gamma_I$ by
permuting the factors in $T^r(H)$, $T^r(s^{N-1} H)$ and the pairs in $I$.
The map
$\gamma_I$, being clearly $\Sigma_r$-equivariant, induces a map
$$\gamma_I : \Gamma^r(s^{N-1}H) \to \Gamma^r H \cdot e_I\,.$$

We define $\phi_I : H^{\otimes n-2r} \to H^{\otimes n}$ by
inserting the element $1\in H^0$ in the position belonging to $I'$
and the map $\psi_I : \Gamma^{\otimes r}(s^{N-1}H) \to H^{\otimes
n}\cdot e_I$ by composition of $\gamma_I$ with the insertion of
the element $1$ in position not in $I_1$.

\medskip
\emph{The central object of this section} is the map
$$\Phi_r : \Gamma^{n-2r} H \otimes \Gamma^r(s^{N-1}H) \to E_1$$
defined by $$\Phi_r (\alpha \otimes \beta) =  \sum_{I\in {\mathcal
P}_r} \phi_I(\alpha)\cdot \psi_I(\beta)\,.$$

\begin{theorem}\label{thm:invariantsk}
 The image of the map
$\Phi = \oplus_{r=0}^{[n/2]} \Phi_r$,
 $$ \Phi : \oplus_{r=0}^{[n/2]}
 \left( \Gamma^{n-2r} H \otimes \Gamma^r(s^{N-1}H)\right) \to E_1$$ is
  the space $E_1^{\Sigma_n}$ of $\Sigma_n$-invariant elements of $E_1$.
\end{theorem}

\begin{proof} Recall
$E_1 \cong \oplus_{G \in {\cF_n}} H^{\otimes l(G)} \cdot G\,.
$
where $l(G)$ denotes the number of components of $G$.

Let ${\mathcal S} = ( s_1, s_2, \ldots , s_q)$ be a sequence of
integers with $1 \leq s_1  \leq s_2 \leq \cdots \leq s_q$. We
denote by $\Gamma_{\mathcal S}$ the set of elements in $\cF_n$ with $q$
components having respectively $s_1, s_2, \ldots ,$ and $ s_q$ vertices.
The direct sum $\oplus_{G\in \Gamma_{\mathcal S}} H^{\otimes q} \cdot
G$ is a $\Sigma_n$-invariant subspace. If $a\in E_1^{\Sigma_n}$,
then we may decompose $a$ in $a= \sum_{\mathcal S} a_{\mathcal S}$,
$a_{\mathcal S}
\in (\oplus _{G\in
\Gamma_{\mathcal S}} H^{\otimes q} \cdot G)^{\Sigma_n}$.

Let fix ${\mathcal S}$ and suppose $s_q >1$. For any sequence $1 \leq
n_1 < n_2 < \cdots <n_{s_q}\leq n$ we write $\Gamma_{\mathcal S} =
\Gamma_1 \cup
 \Gamma_2$ where $\Gamma_1$ denotes the set of
forest in $\Gamma_{\mathcal S}$ in which the vertices $n_1, \ldots ,
n_{s_q}$ are vertices of a same tree. We consider the group
$\Sigma_{s_q}$ consisting of the permutations of the elements
$n_1, \ldots , n_{s_q}$. This is a subgroup of $\Sigma_n$ and the
subspaces $\oplus_{G\in \Gamma_1} H^{\otimes q}\cdot G$ and
$\oplus_{G\in \Gamma_2} H^{\otimes q}\cdot G$ are invariant by the
action of $\Sigma_{s_q}$.
By \thmref{thm:invariantsRN}
and its interpretation in terms of trees, one has
$ (\oplus_{G\in \Gamma_1} G)^{\Sigma_{s_q}} = 0$ if $s_q>2$. Therefore,
by arguing on each element of $H^{\otimes q}$,
one deduces
$(\oplus_{G\in \Gamma_1}
H^{\otimes q}\cdot G)^{\Sigma_{s_q}}=0$ if $s_q>2$.
Since this is true for any sequence
$n_1 < n_2
\cdots < n_{s_q}$, we have
$$(\oplus_{G\in \Gamma_{\mathcal S}}
H^{\otimes q}\cdot G)^{\Sigma_n} = 0\;{\text{and}}\;E_1^{\Sigma_n} \cong
(\oplus_{G \in {\mathcal G}_1} H^{\otimes l(G)}
\cdot G)^{\Sigma_n}\,,$$ where ${\mathcal G}_1$ denotes the set of
forests $G$ in which each component  has at most two vertices.

Let $G\in \cG_1$. Denote by
$I=\{ (i_1,j_1),
(i_2,j_2),\ldots , (i_r,j_r)\}$
the element of ${\mathcal P}_r$ built from the components with exactly
two elements. We have clearly an isomorphism of $\Sigma_n$-vector spaces
$$ \oplus_{G \in {\mathcal G}_1} H^{\otimes l(G)}
\cdot G = \oplus_{r=0}^{[n/2]} \oplus_{I\in {\mathcal P}_r} H^{\otimes
n-2r} \otimes H^{\otimes r} \cdot e_I\,,$$ where the components in
$H^{\otimes n-2r}$ correspond to the components not in $I'$, and
the components in $H^{\otimes r}$ correspond to the components in
$I'$, ordered  by $i_1 < i_2 < \cdots < i_r$.
With the previous notation, we thus have
$$E_1^{\Sigma_n} = \oplus_{r=0}^{[n/2]} \left(\oplus_{I\in {\mathcal
P}_r} H^{\otimes n-2r} \otimes H^r\cdot e_I\right)^{\Sigma_n}\,.$$
If we fix $I \in {\cP}_r$, then each permutation $\sigma$ of
the set $\{1, \cdots , n\} \backslash I'$ and each permutation
$\tau$ of the set of pairs in $I$ preserve $H^{\otimes n-2r}
\otimes H^r \cdot e_I$ and $\oplus_{J\neq I} H^{\otimes
n-2r}\otimes H^{\otimes r} \cdot e_J$. Therefore,
$$\left( \oplus_{I\in {\mathcal P}_r} H^{\otimes n-2r} \otimes
H^r\cdot e_I\right)^{\Sigma_n} \subset \oplus_{I\in {\mathcal P}_r}
\Gamma^{n-2r} H \otimes \Gamma^r H\cdot e_I\,.$$
From the identification
$(E_1,d_1)= \sum_{G\in \cF_n} H^{\otimes l(G)}\cdot G$,
we observe that if $\alpha \otimes \beta\cdot e_I \in \Gamma^{n-2r} H
\otimes
\Gamma^r H\cdot e_I$ and $\sigma\in \Sigma_n$, then $\sigma
(\alpha\otimes \beta \cdot e_I) = \alpha\otimes \beta \cdot
e_{\sigma (I)}$ and $\sigma(I)\in \cP_r$. This implies that $ \Phi:
\oplus_{r=0}^{[n/2]}
 \left( \Gamma^{n-2r} H \otimes \Gamma^r(s^{N-1}H)\right) \to
E_1^{\Sigma_n}$ is an isomorphism.
\end{proof}

\emph{We define now a
 multiplication $\mu $ on
$\oplus_{r=0}^{[n/2]}
 \left( \Gamma^{n-2r} H \otimes \Gamma^r(s^{N-1}H)\right)$
that makes $\Phi$ an isomorphism of algebras.}

\medskip
The graded vector space
 $\Gamma^* H =\oplus_{m^\geq 0} \Gamma^m H$ is a sub-Hopf
 algebra of the Hopf algebra $T(H)$:
\begin{itemize}
\item[---] the multiplication, $*:
\Gamma^pH\otimes \Gamma^qH \to \Gamma^{p+q}H$,  is
 the shuffle product;
\item[---] the comultiplication $\nabla : \Gamma^mH \to
 \oplus_{p+q=m}\Gamma^pH\otimes \Gamma^qH$ is defined by\\
$\nabla = \oplus_{p+q=m}\nabla_{p,q}$, \hspace{1cm}
 $\nabla_{p,q}(x_1x_2\ldots x_m)= x_1\ldots x_p\otimes x_{p+1}
 \ldots x_m$.
\end{itemize}

By multiplying 2 by 2 the elements of $T^{2r}(H)$, we define a
linear map
 $T^{2r}(H) \to T^r(H)$
 $$x_1 \otimes \ldots \otimes x_{2r} \mapsto (x_1\bullet x_2)\otimes
\ldots
 \otimes (x_{2r-1}\bullet x_{2r})$$
 that restricts to a   map
 $$ red : \Gamma^{2r} H\to \Gamma^r H\,.$$

 Finally, we quote two other structures:
\begin{itemize}
\item[---] with the multiplication
 component by component,
  $\nu : T^n(H)\otimes T^n(H)\to T^n(H)$, the spaces $T^n(H)$ and
$\Gamma^n H$ are
 graded
 commutative algebras.
\item[---] a natural action $\bar\nu$ of the
  algebra $\Gamma^n H$   on $\Gamma^n(s^{N-1}H)$ can be defined by\\
 $(a_1\otimes   \ldots \otimes a_n)\cdot
 (s^{N-1}b_1\otimes \ldots \otimes s^{N-1} b_n) =
 \varepsilon\,\,\,
 s^{N-1}(a_1\bullet b_1)\otimes \ldots \otimes s^{N-1}(a_n\bullet b_n)$,\\
 with $\varepsilon = (-1)^{ \sum_{j=2}^n \vert a_j\vert (\vert
 b_1\vert + \cdots + \vert b_{j-1}\vert ) + \sum_{i=1}^n i\vert
 a_i\vert }$.
\end{itemize}

\noindent
We can now state:

\begin{theorem}\label{thm:invariantskalg}
Define a multiplication $\mu$ on
  $ \oplus_{r=0}^{[n/2]}
\Gamma^{n-2r} H \otimes \Gamma^{r}(s^{N-1}H)$  by
$$\mu : \left[\Gamma^{n-2r} H \otimes \Gamma^{r}(s H)\right] \otimes
\left[\Gamma^{n-2s} H\otimes \Gamma^{s}(s H)\right] \to
\Gamma^{n-2r-2s}H\otimes \Gamma^{r+s}(s H),$$
$$\mu (( \alpha \otimes \beta) \otimes ( \gamma \otimes
\delta)) =\sum_{ij} \varepsilon_{ij}\,\,
\nu(\alpha_i\otimes\gamma_j)\otimes
\bar{\nu}(red (\gamma_j')\otimes\beta)
* \bar{\nu}(red (\alpha_i')\otimes\delta)\,,$$
where  $s=s^{N-1}$, $\nabla_{n-2r-2s,2s}(\alpha) =
\sum_i \alpha_i\otimes \alpha_i'$, $\nabla_{n-2r-2s,
2r}(\gamma) =\sum_j \gamma_j\otimes \gamma_j'$
and $\varepsilon_{ij}$ is
the graded sign of the permutation.
With this structure, $\Phi$ becomes an isomorphism of graded algebras.
\end{theorem}

This is a consequence of

\begin{lemma}
The following diagram is commutative
$$\xymatrix{
[\Gamma^{n-2r} H \otimes \Gamma^{r}(s^{N-1}H)]\otimes
 [\Gamma^{n-2s} H\otimes \Gamma^{s}(s^{N-1}H)]
\ar[rr]^-{\Phi_r\otimes \Phi_s}\ar[d]_{\mu}&&
 E_1\otimes E_1 \ar[d]^{\scriptstyle{mult}}\\
\Gamma^{n-2r-2s} H \otimes \Gamma^{r+s}(s^{N-1}H)
\ar[rr]^-{\Phi_{r+s}}&&E_1
}$$
\end{lemma}

\begin{proof}
We denote by $\cdot$ the operations $\nu$ and $\bar{\nu}$. A  simple
computation gives
$$\renewcommand{\arraystretch}{1.6}
\begin{array}{l}
\Phi_r(\alpha \otimes \beta)\cdot \Phi_{s}(\gamma\otimes\delta)
\\
=
mult\left(\sum_{I\in
\cP_r}\phi_I(\alpha)\cdot \psi_J(\beta),
\sum_{I\in
\cP_s}\phi_J(\gamma)\cdot\psi_J(\delta)\right)\\
= \sum_{ij}\sum_{K\in{\mathcal
P}_{r+s}}
\left(
\sum_{I\cup J = K}
\phi_K(\alpha_i)\cdot \phi_K(\gamma_j)\cdot
\psi_I(\beta\cdot red(\gamma_j'))\cdot \psi_J(red(\alpha_i')\cdot
\delta)\right)\\
= \sum_{ij} \sum_{K\in {\mathcal P}_{r+s}} \phi_K(\alpha_i\gamma_j)
\cdot \psi_K((\beta\cdot red(\gamma_j')*(red(\alpha_i')\cdot
\delta))\\
= \Phi_{r+s} (\mu ((\alpha \otimes \beta)\otimes (\gamma\otimes
\delta))).
\end{array}
\renewcommand{\arraystretch}{1}
$$
\end{proof}

\begin{remark} The
cardinality of the set $\cP_r$ is $\displaystyle{\vert {\mathcal P}_r\vert
=
\frac{n!}{(n-2r)!2^rr!}}$\,.

For justifying this formula, denote by $Q_r$ the set of
ordered $r$-uples of disjoint pairs. We have clearly $\vert
Q_r\vert = r! \vert {\mathcal P}_r\vert$. If we fix an element in
$Q_r$, we obtain an element of $Q_{r+1}$ by choosing 2 elements
between the remaining $n-2r$ variables. Therefore,
$$\vert Q_{r+1}\vert = \frac{(n-2r)(n-2r-1)}{2} \vert
Q_r\vert\,,$$ and thus,
$$\vert{\mathcal P}_{r+1}\vert =\frac{(n-2r)(n-2r-1)}{2(r+1)} \vert
{\mathcal P}_r\vert\,.$$ An induction on $r$ based on the last formula
gives the result.
\end{remark}

\section{The differential graded algebra $C_n^H$, $N$ even}
\label{sec:CnH}

The purpose of this section is to use the existence of an isomorphism
between the subspace of invariants $\Gamma^n H$ and the
exterior algebra $\land^nH$ for having a better description
of the differential graded algebra $(E_1,d_1)^{\Sigma_n}$.
This will impose some restriction on the field $\bk$. \emph{For all this
section we suppose that $\bk=\Q$ or $\bk=\F_p$ with $p>n$.}

Let
  $\rho_n : \land^nH \to \Gamma^nH$
 be the symmetrization map defined by
 $$\rho_n(x_1\land x_2\land \cdots \land x_n) =\frac{1}{n!}\sum_{\sigma\in
 \Sigma_n} \varepsilon_\sigma x_{\sigma (1)} \otimes \cdots \otimes
 x_{\sigma (n)}$$
 where   $\varepsilon_{\sigma}$ is the graded sign of the
 permutation
 $$x_1 \ldots x_n \mapsto x_{\sigma (1)} \ldots x_{\sigma
 (n)}\,.$$
The map $\rho_n$ is an
isomorphism which sends
the usual multiplication
$m_{\land}$ of $\land^n H$ on
the shuffle product on $\Gamma^n H$.
We have now to translate
on
$$ C_n^H= \oplus_{r=0}^{[n/2]}
\land^{n-2r} H \otimes \land^{r}(s^{N-1}H)$$
the different structures we used
for constructing the law of
$\oplus_{r=0}^{[n/2]}
\Gamma^{n-2r} H \otimes \Gamma^{r}(s^{N-1}H)$.

\begin{theorem}\label{thm:multiplicationofcnH}
\emph{The
multiplication law} $$
\left[\land^{n-2r} H \otimes \land^{r}(s H)\right]
\otimes \left[\land^{n-2s} H \otimes \land^{s}(s H)\right] \to
\land^{n-2r-2s} H \otimes \land^{r+s}(s H)$$ is defined by
$$\scriptstyle{
(x_1 \land \ldots \land x_{n-2r}\otimes s y_1\land \ldots\land
s y_r)\cdot ( z_1\land \ldots\land  z_{n-2s}\otimes s t_1\land
\ldots\land  s t_s) = \\  \displaystyle{
{\scriptstyle
\frac{1}{(n-2s-2r)!}
 \displaystyle{\sum_{
\renewcommand{\arraystretch}{0.6}\begin{array}{l}
{\scriptscriptstyle \sigma\in \Sigma_{n-2r}}\\
{\scriptscriptstyle \tau\in\Sigma_{n-2s}}\end{array}}}}
\renewcommand{\arraystretch}{1}}
\scriptstyle{ \varepsilon_{\sigma, \tau}\,\,\alpha_{\sigma,
\tau}\otimes \beta_{\sigma, \tau} \land \gamma_{\sigma, \tau}},}
$$
where:
\begin{itemize}
\item[---] $s=s^{N-1}$ and $\varepsilon_{\sigma, \tau}$ denotes the graded
sign of the permutation,
\item[---] $\alpha_{\sigma, \tau}$,
$\beta_{\sigma, \tau}$ and $ \gamma_{\sigma, \tau}$ are defined
by:

\noindent
$\!\!\begin{array}{l}
\alpha_{\sigma, \tau} = (x_{\sigma(1)}\bullet z_{\tau (1)})
\land \ldots \land (x_{\sigma (n-2r-2s)}\bullet z_{\tau
(n-2r-2s)}),\\[.1cm]
 \beta_{\sigma, \tau} =
s (x_{\sigma(n-2r-2s+1)}\bullet x_{\sigma (n-2r-2s+2)}\bullet
t_1)\land \ldots \land s (x_{\sigma (n-2r-1)}\bullet x_{\sigma
(n-2r)}\bullet t_s),\\[.1cm]
\gamma_{\sigma, \tau} = s (z_{\tau (n-2r-2s+1)}\bullet z_{\tau
(n-2r-2s+2)}\bullet y_1)
\land \ldots \land s (z_{\tau (n-2s-1)}\bullet z_{\tau (n-2s)}\bullet
y_r)\,.
\end{array}$
\end{itemize}
\end{theorem}

\medskip
 We first define \emph{a multiplication $\nu': \land^n H\otimes
\land^n H\to \land^n H$ by}
$$\nu' (x_1\land \cdots \land x_n\otimes y_1\land \cdots \land y_n) =
\sum_{\sigma\in\Sigma_n} \varepsilon_\sigma\,\, (x_1\bullet y_{\sigma
(1)})\land \cdots \land (x_n\bullet y_{\sigma
(n)})\,,$$ which makes commutative the following diagram
$$
\xymatrix{
\land^nH\otimes\land^nH
\ar[rr]^-{\rho_n\otimes\rho_n}\ar[d]_{\nu'}&&
\Gamma^n H\otimes \Gamma^n H
\ar[d]^{\nu}\\
\land^nH\ar[rr]_-{\rho_n}&&\Gamma^n H\\
}$$

This \emph{multiplication extends into an action, denoted $\bar\nu'$
of
$\land^n  H$ on $\land^n(s^{N-1}H)$.} This action corresponds, via
$\rho$ to the action $\bar{\nu}$ of $\Gamma^nH$ on $\Gamma^n(s^{N-1}H)$,
i.e. the following diagram trivially commutes
$$
\xymatrix{
\land^nH\otimes \land^n(s^{N-1}H)
\ar[rr]^-{\rho_n\otimes\rho_n}\ar[d]_{\bar{\nu'}}&&
\Gamma^n H\otimes \Gamma^n(s^{N-1}H)
\ar[d]^{\bar{\nu}}\\
\land^n(s^{N-1}H)\ar[rr]_-{\rho_n}&&\Gamma^n(s^{N-1}H)\\
}$$

\emph{A classical diagonal map $\Delta_{p,q}' : \land^{p+q}H\to
\land^pH\otimes
\land^qH$} is defined by
$$\Delta_{p,q}' (x_1\land \cdots \land x_{p+q}) = \sum_{\sigma\in (p,q)Sh}
\varepsilon_\sigma\, (x_{\sigma (1)}\land \cdots \land x_{\sigma (p)})
\otimes (x_{\sigma (p+1)}\land \cdots \land x_{\sigma (p+q)})\,,$$ where
$(p,q)Sh$ denotes the set of $(p,q)$ shuffles of the set $\{1, 2, \cdots ,
p+q\}$. Since we work in free commutative graded algebras, we have
also
$$\Delta_{p,q}' (x_1\land \cdots \land x_{p+q}) = \frac{1}{p!q!}
\sum_{\sigma\in
\Sigma_{p+q}}\varepsilon_\sigma  (x_{\sigma (1)}\land \cdots \land
x_{\sigma (p)}) \otimes (x_{\sigma (p+1)}\land \cdots \land x_{\sigma
(p+q)})\,.$$ The symmetrization map $\rho$ gives a relation between
$\Delta_{p,q}'$ and $\Delta_{p,q}$ in terms of a commutative diagram
$$\xymatrix{
\land^{p+q}H\ar[rr]^-{\rho}\ar[d]_{\Delta_{p,q}'}&&
\Gamma^{p+q} H\ar[d]^{\Delta_{p,q}}\\
\land^pH\otimes \land^qH
\ar[rr]_-{\rho\otimes\rho}&&
\Gamma^p H\otimes \Gamma^q H\\
}$$

Finally \emph{a reduction map
$$red' :\land^{2r}H  \to \land^rH$$
is defined by}
\begin{eqnarray*}
red'(x_1\land \cdots \land x_{2r}) &=&  2^r\, \sum_{I\in {\mathcal
P}_r}
\varepsilon_I\,  (x_{i_1}\bullet x_{j_1})\land  \cdots \land
(x_{i_r}\bullet x_{j_r})\\
& =& \frac{1}{r!}\, \sum_{\sigma\in \Sigma_n}
\varepsilon_\sigma\,(x_{\sigma (1)}\bullet x_{\sigma (2)}) \land \cdots
\land (x_{\sigma (2r-1)}\bullet x_{\sigma (2r)})\,.
\end{eqnarray*}
This reduction map
makes the following diagram commutative
$$\xymatrix{
\land^{2r}H\ar[rr]^-{\rho_{2r}}\ar[d]_{red'}&&
\Gamma^{2r} H \ar[d]^{red}\\
\land^r H \ar[rr]_{\rho_r}&&\Gamma^r H\\
}$$

\begin{proof}[Proof of \thmref{thm:multiplicationofcnH}]

The operations $\nu'$, $\Delta'$ and $red'$ fit together into the
following commutative diagram:
$$\xymatrix@=6pt{
\land^{2n-r}H\otimes\land^r(sH)\otimes \land^{n-2s}H\otimes
\land^s (sH)
\ar[dd]_{q'}\!\!\ar[rr]^-{\rho^{\otimes 4}}&&\!\!
\Gamma^{n-2r}H \otimes\Gamma^r(sH)\otimes
\Gamma^{n-2s}H\otimes\Gamma^s(sH)
\ar[dd]^{q}\\
&&\\
\land^{n-2r-2s}H\otimes \land^{r}(sH)\otimes \land^s(sH)
\ar[rr]^-{\rho^{\otimes 3}} \ar[dd]_{1\otimes m_{\land}}&&
\Gamma^{n-2r-2s}H\otimes \Gamma^{r}(sH)\otimes \Gamma^s(sH)
\ar[dd]^{1\otimes *}\\
&&\\
\land^{n-2r-2s}H\otimes \land^{r+s}(sH)
\ar[rr]_{\rho^{\otimes 2}}&&
\Gamma^{n-2r-2s}H\otimes \Gamma^{r+s}(sH)\\
}$$
where $sH = s^{N-1}H$, $m$ is the usual multiplication in $\land(-)$, $*$
is the shuffle product, and  $q$ and $q'$ are defined by
 $$q =
(\nu\otimes\bar\nu\otimes\bar\nu)\circ T\circ (1\otimes red\otimes
1\otimes 1\otimes red\otimes 1) \circ(\Delta_{n-2r-2s,2s}\otimes
1\otimes \Delta_{n-2r-2s,2r}\otimes 1)$$
 $$q' =
(\nu'\otimes\bar\nu'\otimes\bar\nu')\circ T'\circ (1\otimes
red'\otimes 1\otimes 1\otimes red'\otimes 1)
\circ(\Delta_{n-2r-2s,2s}'\otimes 1\otimes
\Delta_{n-2r-2s,2r}'\otimes 1)$$
 The morphisms $T$ and $T'$ are the permutation maps
$$T : \Gamma^{l}H\otimes \Gamma^sH\otimes \Gamma^r sH\otimes
\Gamma^{l}H\otimes \Gamma^rH\otimes \Gamma^s sH \to
(\Gamma^{l}H)^{\otimes 2}\otimes \Gamma^rH\otimes
\Gamma^r sH\otimes \Gamma^sH\otimes \Gamma^s sH$$
$$T' : \land^{l}H\otimes
\land^sH\otimes \land^r sH\otimes \land^{l}H\otimes
\land^rH\otimes \land^s sH \to (\land^{l}H)^{\otimes 2}\otimes
\land^rH\otimes \land^r sH\otimes \land^sH\otimes
\land^s sH\,,$$ with $l = n-2r-2s$.

The composition on the right is the
multiplication $\mu$ on
 $ \oplus_{r=0}^{[n/2]}
\Gamma^{n-2r} H \otimes \Gamma^{r}(s^{N-1}H)$, therefore the composition
on the left is a multiplication $\mu'$ on $C_n^H$ making $\rho$ a
morphism of algebras.
To determine explicitely $\mu'$, let $$x_1\land  \cdots \land
x_{n-2r}\otimes sy_1\land \cdots \land sy_r\otimes z_1\land \cdots
\land z_{n-2s}\otimes st_1\land \cdots \land st_s$$ be an element of
$\land^{2n-r}H\otimes\land^r(sH)\otimes
\land^{n-2s}H\otimes \land^s (sH) $. By applying
$\Delta_{n-2r-2s,2s}'\otimes 1\otimes \Delta_{n-2r-2s,2r}'\otimes
1$, we obtain
\begin{flushleft}
$\displaystyle{\frac{1}{(n-2r-2s)!(2s)!(n-2r-2s)!(2r)!}}
 \displaystyle{\sum_{\sigma\in \Sigma_{n-2r}, \tau\in \Sigma_{n-2s}}}
$
\end{flushleft}
\begin{flushright}
$x_{\sigma
(1)}\land \cdots \land x_{\sigma (n-2r-2s)}
 \otimes x_{\sigma (n-2r-2s+1)}\land \cdots\land
x_{\sigma (n-2r)}\otimes sy_1\land \cdots \land sy_r$\\
$\otimes z_{\tau (1)}
\land \cdots \land z_{\tau (n-2r-2s)} \otimes z_{\tau (n-2r-2s+1)}\land
\cdots \land z_{\tau (n-2s)}\otimes st_1\land \cdots \land st_s
\renewcommand{\arraystretch}{1}$
\end{flushright}

We apply now the reduction process $1\otimes red'\otimes 1\otimes
1\otimes red'\otimes 1$ and we obtain
\begin{flushleft}
$\displaystyle{\frac{1}{(n-2r-2s)!(s)!(n-2r-2s)!(r)!}
\sum_{\sigma\in \Sigma_{n-2r}, \tau\in \Sigma_{n-2s}}} $\\
\hskip 1cm $ x_{\sigma
(1)}\land \cdots \land x_{\sigma (n-2r-2s)}$\\
\hskip 1cm $
\otimes (x_{\sigma (n-2r-2s+1)}\bullet x_{\sigma (n-2r-2s+2)})\land
\cdots
\land (x_{\sigma (n-2r-1)}\bullet x_{\sigma (n-2r)})$\\
\hskip 1cm $\otimes sy_1\land
\cdots
\land sy_r
\otimes z_{\tau (1)} \land \cdots \land z_{\tau (n-2r-2s)}$\\
\hskip 1cm $ \otimes
(z_{\tau (n-2r-2s+1)}\bullet z_{\tau (n-2r-2s+2)})\land \cdots \land
(z_{\tau (n-2s-1)}\bullet z_{\tau (n-2s)})
\otimes st_1\land \cdots
\land st_s.$
\end{flushleft}

The last step is the composition with
$(\nu'\otimes\bar\nu'\otimes\bar\nu')\circ T'$ which gives
\begin{flushleft}
$\displaystyle{\frac{1}{(n-2r-2s)!}  \sum_{\sigma\in \Sigma_{n-2r},
\tau\in \Sigma_{n-2s}} \varepsilon_{\sigma, \tau}}$\\
\hskip 1cm $ (x_{\sigma
(1)}\bullet z_{\tau (1)}) \land \cdots \land (x_{\sigma (n-2r-2s)}\bullet
z_{\tau (n-2r-2s)})$\\
 \hskip 1cm $\otimes s(y_1\bullet z_{\tau
(n-2r-2s+1)}\bullet z_{\tau (n-2r-2s+2)} )\land \cdots \land s(y_r
\bullet z_{\tau (n -2s-1)}\bullet z_{\tau (n-2s)} ) $\\
\hskip 1cm $ \otimes s(x_{\sigma
(n-2r-2s+1)}\bullet x_{\sigma (n-2r-2s+2)}\bullet t_1)\land \cdots \land
s(x_{\sigma (n -2r-1)}\bullet x_{\sigma (n-2r)}\bullet t_s).
$
\end{flushleft}
This proves the formula for the product. \end{proof}

Before doing the computation of the differential, we need to determine a
system of generators of this algebra:

\begin{proposition}\label{prop:e0}
Denote by $e_0$ the unit of $H$ and by $\cdot$ the multiplication in
$C_n^H$. Then  $\displaystyle\frac{e_0}{n!}$ is the unit of $C_n^H$ and
$$\left(\frac{e_0^{n-2r}}{(n-2r)!}\otimes  sy_1\land \cdots\land
sy_r\right)\cdot \left(\frac{e_0^{n-2}}{(n-2)!} \otimes sa\right)
= \frac{e_0^{n-2r-2}}{(n-2r-2)!} \otimes sy_1\land \cdots \land
sy_r\land sa\,.$$
Moreover, the algebra is generated by the elements of $\land^nH$ and
by the elements $e_0^{n-2}\otimes sa$.
\end{proposition}

\begin{proof} The two first properties are simple computations
from the definition of the mutiplication.
The last one follows from them and from the following relation
$$\begin{array}{l}
\displaystyle{\left(\frac{e_0^{n-2r}}{(n-2r)!} \otimes sy_1\land
\cdots\land  sy_r\right)\cdot \left( \frac{e_0^{2r}}{(2r)!} \land x_1\land
\cdots\land  x_{n-2r}\right) -} \\[.4cm]
\hskip .9cm\left( x_1\land \cdots \land
x_{n-2r}\otimes sy_1\land
\cdots \land sy_r\right) \in \land H\otimes \left(\land^r
sH\right)^{>l}\,,
\end{array}$$
where the elements $x_i$ are of degree $>0$
and $l$
is the degree of
$sy_1\land\cdots
\land sy_r$. A decreasing induction on
the degree of elements in $\land^r sH$ ends the proof.
\end{proof}

We may also observe directly from \thmref{thm:multiplicationofcnH} that
the multiplication law induced on $\land^n H$ is given by
$$
(a_1\land\cdots\land a_n)\cdot (b_1\land\cdots\land b_n) =
\sum_{\sigma\in\Sigma_n} \varepsilon_{\sigma} \;
(a_{\sigma(1)}\bullet b_1)\land \cdots\land (a_{\sigma(n)}\bullet b_n).$$

 We now arrive to the \emph{description of the differential.}
Recall first that  $H =
\oplus_{p=0}^N H^p$, with $H^N = \bk \Omega$.
 The diagonal class $\Delta$ is the element of $H\otimes H$ defined by
 $$\Delta = \sum_i (-1)^{\vert b_i'\vert}b_i \otimes  b_i'$$
 where  the    $\{ b_i\}$ form an homogeneous basis for $H$ and the $\{ b_i'\}
 $ form the dual basis defined by
 $$b_i\bullet b_j' = \delta_{ij} \Omega\,.$$
 By Poincar\'e duality, for each element $h\in H$ we have
$(h\otimes 1)\cdot \Delta = (1\otimes h)\cdot \Delta$.

\begin{definition}\label{def:generald}
We define a differential $d$ on
$\oplus_{r=0}^{[n/2]}
\land^{n-2r} H \otimes \land^{r}(s^{N-1}H)$ by
$$\begin{array}{l}
 d(x_1\land \cdots \land x_{n-2r}\otimes sy_1\land \cdots \land sy_r)
=\\[.4cm]
\displaystyle{\frac{1}{2}
\sum_{i=1}^r\sum_k (-1)^{\vert b_k'\vert + i + \sum_{j=1}^{i-1}
 \vert y_j\vert} x_1\land \cdots \land x_{n-2r} \land
(y_i\bullet b_k)\land b_k'\otimes sy_1\land \cdots\land
 \widehat{sy_i} \land \cdots \land sy_r\,.}
\end{array}$$

\end{definition}

\begin{proof}[Proof of \thmref{thm:invariantQ}]
We are reduced to prove the commutativity of the
following diagram
$$\xymatrix{
\land^{n-2r}H \otimes  \land^r(sH)\ar[rr]^-{\Phi}\ar[d]_{d}&&
E_1\ar[d]^{d_1}\\
\land^{n-2r+2} H \otimes \land^{r-1}(sH)\ar[rr]^-{\Phi}&&E_1\\
}$$
From \propref{prop:e0}, it is enough to check it in the particular
case $r=1$ and $x_1\land \cdots \land
x_{n-2r} = e_0^{n-2}$.
On one hand, we have:
$$d_1\Phi (e_0^{n-2} \otimes sa) = (n-2)! \sum_{i<j} \sum_k
(-1)^{\vert b_k'\vert} 1 \otimes \cdots \otimes (a\bullet b_k)\otimes
\cdots
\otimes b_k'\otimes \cdots \otimes 1\,,$$ where the terms
$(a\bullet b_k)$ and $b_k'$ are respectively in positions $i$ and $j$. On
the other hand, we have:
$$
\renewcommand{\arraystretch}{1.6}
\begin{array}{ll}
\Phi d(e_0^{n-2}\otimes sa) = & \displaystyle\frac{(n-2)!}{2}
\sum_{i<j} \sum_k (-1)^{\vert b_k'\vert} \,\,1 \otimes \cdots
\otimes
(a\bullet b_k)\otimes \cdots \otimes b_k'\otimes \cdots \otimes 1 \\
&+ \sum_{i<j} \sum_k (-1)^{\vert b_k'\vert + \vert b_k'\vert
\cdot \vert ab_k\vert} \,\,1 \otimes \cdots \otimes  b_k'\otimes
\cdots \otimes (a\bullet b_k) \otimes \cdots \otimes 1\,.
\end{array}
\renewcommand{\arraystretch}{1}
$$
In
the first sum the elements $(a\bullet b_k)$ and $b_k'$ are in positions $i$
and $j$, and in the second sum the elements $b_k'$ and $(a\bullet b_k)$
are in positions $i$ and $j$.

Since $N$ is even and $\Delta = \sum_k (-1)^{\vert b_k'\vert}
b_k\otimes b_k' = \sum_k (-1)^{\vert b_k\vert + \vert b_k'\vert +
\vert b_k\vert } b_k' \otimes b_k$, we have
$$ \sum_k (-1)^{\vert b_k'\vert + \vert ab_k\vert \cdot \vert
b_k'\vert } b_k'\otimes (a\bullet b_k) = \sum_k(-1)^{N + \vert b_k'\vert}
(a\bullet b_k) \otimes b_k' = \sum_k (-1)^{\vert b_k'\vert} (a\bullet
b_k)\otimes b_k'\,.$$ Therefore $\Phi\circ d = d_1\circ \Phi$.
\end{proof}

\begin{proof}[Proof of \thmref{thm:collapseQ}, $N$ even]
Observe that the underlying complex of the cdga $(C_n^H,d)$ we have
described before, is isomorphic to the complex defined by Y.~F\'elix
and {J.-C.~Thomas}
\cite{FT} for the determination of the homology vector space
$H_*(C_n(M);\Q)$.
\end{proof}
\section{The space of unordered configurations in $P^2(\mathbb
C)$}\label{sec:configCP2}

In this section, we apply the previous results to the unordered
configuration space of
$P^2(\C)$ and prove \thmref{thm:configCP2}.

\smallskip
Let
$H = H^*(P^2(\mathbb C);\mathbb Q)$. As quoted in the introduction
(see also \defref{def:generald}), the
direct sum
$(\oplus_n C_n^H,d)$ is, as a complex,  the commutative differential
graded algebra $(\land (x_0,x_1,x_2, y_0, y_1, y_2),d)$ with
$|x_0|=0$, $|x_1|=2$, $|x_2|=4$, $|y_0|=3$, $|y_1|=5$, $|y_2|=7$ and
the differential
$d$ defined by
$$d(x_0)=d(x_1)=d(x_2)=0\,; \; d(y_0) = x_0x_2 + \frac{1}{2} x_1^2\,;\;
d(y_1) = x_1x_2\,; \;d(y_2) = \frac{1}{2}x_2^2\,.$$

\begin{lemma}
A basis of the reduced
cohomology of $(\land (x_0,x_1,x_2, y_0, y_1, y_2),d)$ is given by
  the classes of the cocycles
\begin{enumerate}
\item $x_0^n$, $x_0^{n-1}x_1$ and $x_0^{n-1}x_2$, for $n\geq 1$,
\item $x_0^{n-3}x_1y_1 -2x_0^{n-3}x_2y_0 + 4 x_0^{n-2} y_2$ and
$x_0^{n-3}x_2y_1   -2 x_0^{n-3}x_1y_2$ for $n\geq 3$,
 \item
$x_0^{n-4}x_1x_2y_1 - 2 x_0^{n-4} x_1^2y_2$ for $n\geq 4$.
\end{enumerate}
\end{lemma}

\begin{proof} We consider the relative
minimal model (see \cite{FHT})
$$
\xymatrix@1{
(\land (x_0,0) \ar[r]&  (\land (x_0,x_1,x_2, y_0, y_1, y_2),d)
\ar[r]^-{p}&
(\land ( x_1,x_2, y_0, y_1, y_2),\bar d)\,.
}$$
Since $x_1^2$ and $x_2^2$ form a regular sequence in $\land
(x_1,x_2)$, the canonical projection
$$(\land ( x_1,x_2, y_0, y_1, y_2),\bar d) \to (\land ( x_1,x_2,   y_1 )/(x_1^2, x_2^2),\bar d(y_1) =
x_1x_2)$$ is a quasi-isomorphism. The algebra $\land ( x_1,x_2,
y_1 )/(x_1^2, x_2^2)$ is finite dimensional and a basis of the
reduced cohomology is given by the classes of the cocyles  $x_1$,
$x_2$, $y_1x_1$, $y_1x_2$, $y_1x_1x_2$. The morphism $H^*(p)$ is
surjective because $H^*(p) ( [x_1y_1 -2x_2y_0 + 4 x_0 y_2]) =
[x_1y_1]$, $H^*(p) ([x_2y_1   -2 x_1y_2]) = [x_1y_2]$ and
$H^*(p)([ x_1x_2y_1 - 2   x_1^2y_2]) = [y_1x_1x_2]$. Therefore the
spectral sequence obtained by filtering the complex
$(\land (x_0,x_1,x_2, y_0, y_1, y_2),d)$
by the powers
of the ideal generated by $x_0$ collapses at the $E_2$-term. This
gives the above basis for the cohomology.
\end{proof}

\begin{proof}[Proof of \thmref{thm:configCP2}]
Let fix $n$. Then a linear
basis of the cohomology of
$C_n^H$ is given by the classes
$$
 \renewcommand{\arraystretch}{1.6}
 \left\{
\begin{array}{l}
 x_0^n\,, \; x_0^{n-1}x_1 \,, \;x_0^{n-1}x_2 \,,\hspace{5mm} (n\geq 1)\,,\\
 z_0= x_0^{n-3}x_1y_1 -2x_0^{n-3}x_2y_0 + 4 x_0^{n-2} y_2\,,
\hspace{5mm} (n\geq
 3)\,,\\
 z_1 = x_0^{n-3}x_2y_1   -2 x_0^{n-3}x_1y_2\,, \hspace{5mm} (n\geq
 3)\,,\\
 z_2 = x_0^{n-4}x_1x_2y_1 - 2 x_0^{n-4} x_1^2y_2\,, \hspace{5mm} (n\geq 4)\,.
 \end{array}
 \right.
  \renewcommand{\arraystretch}{1}
 $$

 Denote $t_0 = \displaystyle\frac{x_0^n}{n!}$, $t_1 =
 \displaystyle\frac{x_0^{n-1}x_1}{(n-1)!}$ and $t_2 =
 \displaystyle\frac{x_0^{n-1}x_2}{(n-1)!}$. Then $t_0$ is the
 unit, and $t_1^2 = (3-2n)t_2$.
 Moreover,   $$[t_1]\cdot [z_1] = (n-3)[z_2]$$ for $n\geq 4$ and
 $$[t_1]\cdot [z_0] = (3-2n) [z_1]$$ for $n\geq 3$.
 This gives the statement.
\end{proof}

 \section{Unordered configuration space of odd dimensional
manifolds}\label{sec:oddcase}

\begin{proof}[Proof of \thmref{thm:casimpair}]
The manifold $M$ being odd dimensional, one has $e_{ij}=-e_{ji}$ and the
elements $e_{ij}$ are of even degree. We have also to take in account the
relation $e^2_{ij}=0$ which is automatic for degree reason when $M$ is
even dimensional.

Let
$\alpha\in E_1^{\Sigma_n}= \left(H^{\otimes n}\otimes \land
e_{ij}/I\right)^{\Sigma_n}$. We decompose $\alpha$ in
$\alpha=\alpha_{n-1}+\cdots \alpha_0$ where
$\alpha_r\in K_r = H^{\otimes n}\otimes e_{rn}\otimes \land (e_{ij})/I$
with the remaining $e_{ij}$ such that
$(i,j)\not\in\{
 (r,n), (r+1,n)\cdots ,(n-1,n)\}$. For instance,
$\alpha_0\in H^{\otimes n}\otimes \land (e_{ij})_{i<j<n}/I$.

Let $s$ be the greatest integer such that $\alpha_s\neq 0$. The image
$\tau(\alpha)$ of
$\alpha$ by the permutation $(s,n)$ satisfies
$\tau(\alpha)= - \alpha_s + \beta$ with $\beta \in \oplus_{i\neq s}K_i$.
Therefore, if $\alpha\in E_1^{\Sigma_n}$, we have $\alpha_s=0$ and, by
induction,
$\alpha\in (H^{\otimes n})^{\Sigma_n}= \Gamma^n H$.

The first term of the spectral sequence $\cE(M,n,\bk)^{\Sigma_n}$
is \emph{concentrated in filtration degree~0}. Therefore, the spectral
sequence collapses and $\Gamma^n H$ (which coincides with $\land^n H$
with our hypothesis on the field $\bk$) is the algebra of cohomology
$H^*(C_n(M);\bk)$.
\end{proof}

As a graded vector space this corresponds to the
result previously obtained by C.F.~B\"o\-di\-gheimer, F.~Cohen and
L.~Taylor (\cite{BCT}). Here we prove that, in the case of an odd
dimensional manifold $M$, the cohomology
\emph{algebra} of unordered configurations depends only on Betti number of
the manifold $M$.


\end{document}